\documentclass[aop]{imsart}

\usepackage{bm}
\usepackage{amsmath, latexsym, amsfonts, amssymb, amsthm}
\usepackage{graphicx, color,hyperref,dsfont,epsfig,caption,wrapfig,subfig}
\usepackage{pstricks}
\usepackage{datetime}
\usepackage{movie15}
\hypersetup{
    linktoc=page,
    linkcolor=red,          
    citecolor=blue,        
    filecolor=blue,      
     urlcolor=cyan,
    colorlinks=true           
}

\numberwithin{equation}{section}

\newtheorem{theorem}{Theorem}[section]
\newtheorem{lemma}[theorem]{Lemma}

\newtheorem{corollary}[theorem]{Corollary}

\newtheorem{remark}[theorem]{Remark}
\newtheorem{definition}[theorem]{Definition}

\theoremstyle{definition}

\definecolor{remi}{rgb}{0,0,0}

\renewcommand{\tilde}{\widetilde}          
\DeclareMathSymbol{\leqslant}{\mathalpha}{AMSa}{"36} 
\DeclareMathSymbol{\geqslant}{\mathalpha}{AMSa}{"3E} 
\DeclareMathSymbol{\eset}{\mathalpha}{AMSb}{"3F}     
\renewcommand{\leq}{\;\leqslant\;}                   
\renewcommand{\geq}{\;\geqslant\;}                   

\newcommand{\R}{\mathbb{R}}

\newcommand{\E}{\mathds{E}}
\renewcommand{\P}{\mathds{P}}

\newcommand{\caB}{{\mathcal B}}

\newcommand{\C}{\mathbb{C}}

\renewcommand{\P}{\mathds{P}}

\usepackage{xparse}

\def\bi{\begin{itemize}}
\def\ei{\end{itemize}}

\def\bnum{\begin{enumerate}}
\def\enum{\end{enumerate}}
\def\<#1{\langle #1 \rangle}

\begin{document}

\begin{frontmatter}

\title{The Tail expansion of Gaussian multiplicative chaos and the Liouville reflection coefficient}
\runtitle{Tail of GMC measures}

\begin{aug}
  \author{\fnms{R\'emi  }  \snm{Rhodes}\corref{}\thanksref{t2}\ead[label=e1]{remi.rhodes@univ-amu.fr}}
  \and
  \author{\fnms{Vincent}  \snm{Vargas}\thanksref{t2}%
  \ead[label=e3]{Vincent.Vargas@ens.fr}%
  \ead[label=u1,url]{http://www.foo.com}}

  \thankstext{t2}{Research supported in part by ANR grant Liouville (ANR-15-CE40-0013).}

  \runauthor{R. Rhodes and V. Vargas}

  \affiliation{Universit\'e Aix-Marseille and \'Ecole Normale Sup\'erieure de Paris}

  \address{Universit\'e Aix-Marseille
Bureau R326,\\
Institut de MathŽmatiques de Marseille (I2M),\\
Centre de MathŽmatiques et Informatique (CMI)\\
Technop\^ole Ch\^ateau-Gombert\\
39, rue F. Joliot Curie
13453 Marseille Cedex 13,\\ 
          \printead{e1}}

  \address{ENS Ulm, DMA, \\45 rue d'Ulm, \\ 75005 Paris, France,\\
          \printead{e3}}

\end{aug}

\begin{abstract}
 In this short note, we derive a precise tail expansion for  Gaussian multiplicative chaos (GMC) associated to the 2d GFF on the unit disk with zero average on the unit circle (and variants). More specifically, we  show that to first order the tail is a constant times an inverse power with an explicit value for the tail exponent as well as an explicit value for the constant in front of the inverse power; we also provide a second order bound for the tail expansion. The main interest of our work consists of two points. First, our derivation is based on a simple method which we believe is universal in the sense that it can be generalized to all dimensions and to all log-correlated fields. Second, in the 2d case we consider, the value of the constant in front of the inverse power is (up to explicit terms) nothing but the Liouville reflection coefficient taken at a special value. The explicit computation of the constant was performed in the recent rigorous derivation with A. Kupiainen of the DOZZ formula \cite{KRV1,KRV}; to our knowledge, it is the first time one derives rigorously an explicit value for such a constant in the tail expansion of a GMC measure.  We have   deliberately kept this paper short to emphasize the method  so that it becomes an easily accessible toolbox  for computing tails in GMC theory.
\end{abstract}

\begin{keyword}[class=MSC]
\kwd[Primary ]{60G57}
\end{keyword}

\begin{keyword}
\kwd{Gaussian multiplicative chaos}\kwd{tail expansion}\kwd{reflection coefficient}\kwd{DOZZ formula}
\end{keyword}

\end{frontmatter}


%
%
%
%
%



\maketitle
\tableofcontents

\section{Introduction}

Gaussian multiplicative chaos (GMC) measures are widespread in probability and statistical physics; indeed they appear in a wide variety of contexts and in particular in the fields of finance, number theory, Liouville quantum gravity and turbulence (see \cite{review} for references). In view of the broad applications of GMC theory, it is very natural to study these measures in great detail. The foundations of GMC theory were laid in Kahane's 1985 seminal work \cite{cf:Kah}. If $\Omega$ is some open subset of $\R^d$ then the theory of GMC enables one to define random measures of the form
\begin{equation}\label{defchaos1}
M_\gamma(dx)= e^{\gamma X(x)-\frac{\gamma^2}{2} \E[X(x)^2]} dx
\end{equation}
where $dx$ denotes the Lebesgue measure and $X$ is a log-correlated Gaussian field, i.e. a centered Gaussian field with covariance 
\begin{equation}\label{defkernel}
\E[X(x)X(y)]= \ln \frac{1}{|x-y|}+ f(x,y),\quad x,y\in \Omega
\end{equation}
with $f$ some smooth and bounded function. Of course, definition \eqref{defchaos1} is only formal since $X$ is not defined pointwise hence the measure $M_\gamma$ can only be defined via a regularization procedure. More specifically, the measure  $M_\gamma$ is defined via the limit in probability of the sequence 
\begin{equation}\label{defchaos}
M_{\gamma,\epsilon}(dx)= e^{\gamma X_\epsilon(x)-\frac{\gamma^2}{2} \E[X_\epsilon(x)^2]} dx
\end{equation}
where $(X_\epsilon)_{\epsilon>0}$ is a reasonable family of smooth Gaussian fields converging towards $X$ when $\epsilon$ goes to $0$: see Berestycki's very simple and elegant approach \cite{Ber} for an introduction to GMC and an account on the above issues of convergence.

Kahane \cite{cf:Kah} proved in 1985 that the measure $M_\gamma$ is different from $0$ if and only if $\gamma^2 <2d$\footnote{We will only consider this case in this note though there has been much progress recently in understanding the critical case $\gamma^2=2d$.}; moreover, a standard result in GMC theory  is the following condition on existence of moments (see \cite[section 2]{review}): if $\mathcal{O}$ is some nonempty and bounded open subset then 
\begin{equation}\label{eq:exismom}
 \E[  M_\gamma(  \mathcal{O} )^p ]< \infty \quad   \Longleftrightarrow \quad p< \frac{2d}{\gamma^2}.
 \end{equation}
  In view of  \eqref{eq:exismom}, it is natural to seek the exact tail behaviour of $M_\gamma(  \mathcal{O}) $. This was achieved in the beautiful work by Barral and Jin   \cite{BaJin} in the 1d case and for a covariance kernel  of the form $\E[X(x)X(y)]=\ln\frac{1}{|x-y|}$ on the interval $[0,1]$. More specifically, the work of Barral and Jin established that
 \begin{equation}\label{tailM}
 \P(  M_\gamma [0,1] > t ) = \frac{C_\star}{t^{\frac{2}{\gamma^2}}} +o(\frac{1}{t^{\frac{2}{\gamma^2}}})
 \end{equation}
 where the constant $C_\star$ is  given by 
 \begin{equation}\label{defconstant}
 C_\star= \frac{2\gamma^2}{2-\gamma^2}  \frac{\E[  M_\gamma [0,1]^{\frac{2}{\gamma^2}-1} M_\gamma [0,\frac{1}{2}]- M_\gamma [0,\frac{1}{2}]^{\frac{2}{\gamma^2}}       ]}{ \ln 2}.
 \end{equation}
 Unfortunately, it is not obvious how to generalize the work \cite{BaJin} to higher dimensions and other kernels of the type \eqref{defkernel}  because their argument is based on a functional relation, which is obtained from a specific geometric  representation (called  cone construction)  of the Gaussian field with kernel  $\E[X(x)X(y)]=\ln\frac{1}{|x-y|}$ in dimension 1.  Moreover, the approach of Barral and Jin does not provide an explicit value for the constant $C_\star$.   Let us also mention the appendix of Lebl\'e-Serfaty-Zeitouni  \cite{LSZ} in a slightly different framework, which contains a description of the tail of the modulus of complex GMC (for $\beta\in (0,\sqrt{2})$)
 $$M_\beta(D):=\lim_{\epsilon\to 0}\epsilon^{-\frac{\beta^2}{2}}\int_De^{i\beta X_\epsilon(x)} dx$$ with $X$ a full plane Gaussian Free Field in dimension 2 and $D$ a bounded open subset of the plane. This type of GMC possesses moments of all orders, which are related to  Coulomb gas computations. This makes possible to compute the tail with a  value of the constant  given by   a variational formula (whose solution is not explicit).  
 
Let $\mathcal{O}$ denote an open subset of $\R^2$ with a $C^1$ boundary\footnote{Our techniques could perhaps handle more general cases but for the sake of simplicity we consider only the case of a smooth boundary: see Remark \ref{C1assump}.}.  The aim of this work is to introduce a simple method  to compute the tail of $M_\gamma (\mathcal{O})$  together with an exact value for the constant in the leading order term. We will consider the case of  the unit disk  $D$ equipped with  a Gaussian Free Field (GFF) $X$ with vanishing mean over the unit circle  (eventually augmented by an independent Gaussian perturbation; this covers for instance the important case of the GFF with Neumann boundary condition). Indeed, though we believe our method can be generalized to all dimensions and all kernels, we stick to the 2d GFF setting to keep this note rather short. Let us stress furthermore that these generalizations may eventually raise serious additional technicalities.

We further mention that some material in this note is inspired by the work of Duplantier-Miller-Sheffield \cite{DMS} and also by the tail estimates which appear in our recent proof with Kupiainen of the DOZZ formula \cite{KRV1,KRV}. In fact, our main result (Theorem \ref{th1} below) can be seen as a strengthening of the convergence results which underly the construction of the so-called quantum sphere in \cite{DMS}. However, this note is mostly self-contained: it requires no a priori knowledge of the  paper \cite{DMS} and we recall the definition of the reflection coefficient of Liouville conformal field theory (LCFT hereafter)   before stating the main result.

\subsection{The Liouville reflection coefficient}
 
 From now on we restrict to the case of dimension $2$. Consider $\gamma \in (0,2)$ and define $Q=\frac{\gamma}{2}+\frac{2}{\gamma}$.
In this subsection, we introduce the (unit volume) reflection coefficient of LCFT.   It was defined in \cite{KRV} where it plays an important role in the proof of the DOZZ formula.  This coefficient shows up when one analyzes the large values of GMC measures with a singularity, i.e. when one looks at variables of the type
$$\int_{B(v,r)}\frac{1}{|x-v|^{\gamma\alpha}}M_\gamma(d^2x).$$
For this integral to exist, the exponent must satisfy $\alpha<Q$ \cite{DKRV}. It turns out that when $\alpha$ is large enough (in fact $\alpha>\frac{\gamma}{2}$), large values of this random variable mostly come from the singularity at the point $v$ and the reflection coefficient quantifies this statement (more details later).  Therefore it can intuitively be understood as a coefficient   of mass localization in GMC theory. Our proof for the tail is based on a localization trick (see subsection \ref{sub:localization}), which allows us to express the tail of GMC in terms of singular integrals as above. This is the reason why the reflection coefficient is  instrumental in identifying the constant in the tail of GMC measures.

In order to introduce the reflection coefficient, we first recall basic material introduced in \cite{zipper}.  For all $\alpha<Q$, we define the process $ \mathcal{B}^\alpha_s$ 
\begin{equation}\label{BMneg}
 \mathcal{B}^\alpha_s = \left\{
 \begin{array}{ll}
  B^\alpha_{-s} & \text{if } s < 0\\
    \bar{B}^\alpha_{s} & \text{if } s >0 \end{array} \right.
\end{equation}
where $(B^{\alpha}_s)_{s \geq 0},(\bar B^{\alpha}_s)_{s \geq 0}$ are  two independent standard Brownian motions with negative drift $\alpha-Q$ and conditioned to stay negative. We also consider an independent centered Gaussian field $Y$ defined for $s\in\R$ and $\theta\inÊ[0,2\pi]$ with  covariance 
\begin{equation}\label{covlateral}
\E[  Y(s,\theta) Y(s',\theta') ]  = \ln \frac{e^{-s}\vee e^{-s'}}{|e^{-s}e^{i \theta} - e^{-s'} e^{i \theta'} |}
\end{equation}
and associated GMC measure
\begin{equation*}
N_\gamma (ds d \theta)=e^{\gamma Y(s,\theta)-\frac{\gamma^2}{2} \E[Y(s,\theta)^2]} ds d\theta.
\end{equation*}
We introduce the integrated chaos measure with respect to $Y$ 
\begin{equation}\label{DefZfirst}
Z_s = \int_0^{2 \pi} e^{\gamma Y(s,\theta)-\frac{\gamma^2}{2} \E[Y(s,\theta)^2]} d\theta.  
\end{equation}
This is a slight abuse of notation since the process $Z_s$ is not a function (for $\gamma \geq \sqrt{2}$) but rather a generalized function; in this setup, $Z_sds$ is a stationary 1d random measure,  meaning that this random measure has same law as its   pushforward by any translation. 

Now, we define the (unit volume) reflection coefficient $\bar{R}(\alpha)$ for all $\alpha \in (\frac{\gamma}{2},Q)$ by the following formula
\begin{equation}\label{firstdefunitR}
 \bar{R}(\alpha)=  \E\Big[ \left ( \int_{-\infty}^\infty  e^{   \gamma \mathcal{B}_s^\alpha  } Z_s ds  \right )^{\frac{2}{\gamma}  (Q-\alpha)}\Big].
\end{equation}

Notice that the condition $\alpha \in (\frac{\gamma}{2},Q)$ ensures that $\frac{2}{\gamma}  (Q-\alpha) < \frac{4}{\gamma^2}$ hence one can show that  $\bar{R}(\alpha)$ is well defined for all  $\alpha \in (\frac{\gamma}{2},Q)$. Indeed, for all $\alpha<Q$ and $p \in (0,\frac{4}{\gamma^2})$ the following holds (see \cite[Lemma 2.8]{KRV})
\begin{equation}\label{momentcondR}
 \E\Big[ \left ( \int_{-\infty}^\infty  e^{   \gamma \mathcal{B}_s^\alpha  } Z_s ds  \right )^{p}\Big] < \infty.
\end{equation}

Finally, one of the main results of \cite{KRV} is an integrability result for $\bar{R}$. Indeed, one has the following remarkable explicit formula for $\bar{R}$ 
\begin{equation}\label{refformula}
\bar{R}(\alpha)= -\frac{(\pi l(\frac{\gamma^2}{4}))^{\frac{2}{\gamma}(Q-\alpha)}}{\frac{2}{\gamma} (Q-\alpha)}   \frac{\Gamma (  - \frac{\gamma}{2}(Q-\alpha)  )}{\Gamma (   \frac{\gamma}{2}(Q-\alpha)  ) \Gamma (  \frac{2}{\gamma}(Q-\alpha)  )}
\end{equation}
where $\Gamma$ is the standard Gamma function and $l(x)=\Gamma(x)/\Gamma(1-x)$.
\vspace{0.3 cm}

 Now that we have introduced $\bar{R}$, we can announce the organization of the paper. First, we state our main result in section \ref{setup} and then proceed with the proof in section \ref{proof}. The final section of the paper discusses how the method can be extended to more general situations.
 
\section{Setup  and Main Results} \label{setup}
 We consider  the open unit disk  $D\subset\C$ equipped with  a   full plane  Gaussian Free Field (GFF) $X$   with vanishing mean over the unit circle. This is a centered Gaussian field with covariance
 \begin{equation}\label{hatGformulafirst}
\forall x,y\in D,\quad \E [ X(x)X(y)] =\ln\frac{1}{|x-y|}.
 \end{equation}
\begin{remark}
Recall that the full plane GFF is a centered Gaussian distribution $\tilde{X}$ defined up to additive constant with covariance $\ln\frac{1}{|x-y|}$. Thus it makes sense to consider $X:=\tilde{X}-\frac{1}{2\pi i}\oint_{|z|=1}\tilde{X}(z)\frac{dz}{z}$, for which the covariance is easily identified to be given by 
$$\E[X(x)X(y)]=\ln\frac{1}{|x-y|}+\mathbf{1}_{\{|x|\geq 1\}}\ln|x|+\mathbf{1}_{\{|y|\geq 1\}}\ln|y|.$$
When restricted to $D$, this gives \eqref{hatGformulafirst}.
\end{remark}
 
For $\gamma\in (0,2)$, we consider the GMC measure
\begin{equation}\label{GMC}
M_\gamma(d^2z):= e^{\gamma X(z)-\frac{\gamma^2}{2}\E[X(z)^2]} d^2z
\end{equation}
where $d^2z$ is the standard Lebesgue measure in the plane (recall that \eqref{GMC} is defined as the limit of \eqref{defchaos} when $\epsilon$ goes to $0$). In the sequel, we adopt the following standard convention: for $\beta \in \R$, we denote $o(t^{\beta})$ (respectively $O(t^{\beta})$) a quantity of the form $\epsilon_t t^{\beta}$ where $\epsilon_t$ goes to $0$ (resp. remains bounded) as $t$ goes to infinity. Denote by $|B|$ the Lebesgue measure of the Borel measurable set $B$ and by $\bar{A}$ the closure of the set $A$.  

Our main result is: 
 \begin{theorem}\label{th1}
 Set  
\begin{equation}\label{defp0}
p_0:=\sqrt{\frac{(2-\tfrac{\gamma^2}{2})^2}{\gamma^4}+\frac{2}{\gamma^2}}+\frac{(2-\tfrac{\gamma^2}{2})}{\gamma^2}.
\end{equation}
For any $\delta\in (0,\frac{1+p_0-\frac{4}{\gamma^2}}{2+p_0})$ and  for any open set $\mathcal{O}\subset D$ with a $C^1$ boundary, we have 
\begin{equation}\label{mainresultpap}
\P(M_\gamma(\mathcal{O})>t)= \frac{\frac{2}{\gamma}(Q-\gamma)}{\frac{2}{\gamma}(Q-\gamma)+1}   \frac{\bar{R}(\gamma)}{t^{\frac{4}{\gamma^2}}} |\mathcal{O}|+o(t^{-\frac{4}{\gamma^2}-\delta})
\end{equation}
where $\bar{R}(\gamma)$ is the (unit volume) reflection coefficient of LCFT evaluated at $\gamma$. It has explicit expression
\begin{equation}\label{reflexDOZZ}
\bar{R}(\gamma)= -\frac{(\pi l(\frac{\gamma^2}{4}))^{\frac{2}{\gamma}(Q-\gamma)}}{\frac{2}{\gamma} (Q-\gamma)}   \frac{\Gamma (  - \frac{\gamma}{2}(Q-\gamma)  )}{\Gamma (   \frac{\gamma}{2}(Q-\gamma)  ) \Gamma (  \frac{2}{\gamma}(Q-\gamma)  )}\quad\text{ with }\quad Q=\frac{2}{\gamma}+\frac{\gamma}{2}
\end{equation}
and the function $l$ is given by a ratio of Gamma functions $l(x)=\Gamma(x)/\Gamma(1-x)$.
 \end{theorem}
The origin of the value $p_0$ will be more transparent in Remark \ref{spectrum} below. Also, notice that $1+p_0-\frac{4}{\gamma^2}>0$ so that our statement is not trivially empty.
\begin{remark}
In a previous version of this note, we stated \eqref{mainresultpap} with a better bound on $\delta$ than the present condition; unfortunately this was due to a slight mistake in the proof. Also, let us mention that we do not know what is the optimal bound on $\delta$ such that \eqref{mainresultpap} holds.
\end{remark}

More than the result, our method of proof makes us believe in a higher level  of generality for this result. For general GMC measures in $2d$ with respect to a Gaussian field with covariance \eqref{defkernel}: we expect that the $\ln\frac{1}{|x-y|}$ part is responsible for the appearance of  the reflection coefficient term whereas the perturbation term $f(x,y)$ in \eqref{defkernel} produces an integral which depends only on the on-diagonal values of $f$.

As an illustration of this discussion we claim:
 \begin{corollary}\label{th2}
Assume the GMC measure \eqref{GMC} is defined with respect to the Gaussian distribution $X$ with covariance \eqref{defkernel} for some   function $f$  that is locally H\"{o}lder and positive definite.  For any open set $\mathcal{O}\subset \bar{\mathcal{O}}\subset D$ with a $C^1$ boundary, we have 
\begin{equation}\label{mainresultpap2}
\P(M_\gamma(\mathcal{O})>t)=  \left (  \int_\mathcal{O}  e^{ \frac{4}{\gamma}(Q-\gamma) f(v,v) } d^2v \right ) \frac{\frac{2}{\gamma}(Q-\gamma)}{\frac{2}{\gamma}(Q-\gamma)+1}   \frac{\bar{R}(\gamma)}{t^{\frac{4}{\gamma^2}}}  +o(t^{-\frac{4}{\gamma^2} })
\end{equation}
where $\bar{R}(\gamma)$ is still the (unit volume) reflection coefficient of LCFT evaluated at $\gamma$. 
 \end{corollary}

\begin{remark}
One can apply Corollary \ref{th2} with $f(x,y)=\ln\frac{1}{|1-x\bar{y}|}$ to get the right tail of GMC based on the GFF $X$ with Neumann boundary condition in $D$, which has covariance $\E[X(x)X(y)]=\ln\frac{1}{|x-y||1-x\bar{y}|}$. Indeed  $f(x,y)=\ln\frac{1}{|1-x\bar{y}|}$ is positive definite\footnote{The fact that $f$ is positive definite can be seen after a series expansion of the $\ln$.} and locally H\"older in $D$.
\end{remark}

Actually we even believe that the above structure of the tail of GMC measures  is universal in the sense that in all dimensions  there should be  an analogue of  the unit  volume  reflection coefficient, with a probabilistic representation similar to  \eqref{firstdefunitR}, such that the tail is given by \eqref{mainresultpap2} (see section \ref{beyond} for further discussion). A specific feature of the 2d case is that there exists an explicit analytic expression \eqref{reflexDOZZ} established in \cite{KRV1,KRV} for the expectation formula  \eqref{firstdefunitR}.

Let us comment on the physics literature. To our knowledge, an explicit tail expansion for GMC was derived (at the level of physics rigor) in the papers \cite{fybu, FLeR} in the 1d case for two specific log-correlated models: the circular case in \cite{fybu}  and the unit interval case in \cite {FLeR} (with exact logarithmic correlations, i.e. the case where $f=0$ in \eqref{defkernel}). Their derivation is based on exact integrability results for GMC hence their results are much stronger than just tail expansions for these specific GMC measures; however, these works do not adress the derivation of tail expansions for 1d GMC associated to general logarithmic kernels or more importantly for GMC in higher dimensions. We refer to the section \ref{beyond} below for more on this.

Finally, let us mention that the  our main result  Theorem \ref{th1}  (together with its Corollary \ref{th2}) reinforces/strengthens the convergence results used to define the unit volume quantum sphere \cite[definition 4.21]{DMS} stated by Duplantier-Miller-Sheffield. Indeed the definition of quantum sphere in \cite{DMS} is the following: consider the random field $h$ defined on the cylinder $\R \times [0,2\pi]$ by the sum of two independent processes $h_1(s)+Y(s,\theta)$ such that
\begin{itemize}
\item
The  radial part  $h_1$ is  $\frac{2}{\gamma} \ln U $ where $U$ is chosen according to the a Bessel excursion measure with index $\delta= 4-\frac{8}{\gamma^2}$ (and reparametrized to have quadratic variation $ds$).
\item
The non radial part $Y$ has the law given by \eqref{covlateral}. 
\end{itemize}
The (unit volume or area) quantum sphere is obtained by conditioning the measure $e^{\gamma h} ds d\theta $ to have volume $1$. By definition,   $\frac{2}{\gamma} \ln U $ is distributed according to $\mathcal{B}^\gamma+ \frac{2}{\gamma} \ln e^{\star}$ where $(e^{\star})^2$ (which is the maximum of the excursion) is distributed according to the (infinite) measure $  u^{\frac{\delta}{2}-2} du$ (recall that $\delta= 4-\frac{8}{\gamma^2}$). Recalling the definition \eqref{DefZfirst} of $Z_s$, we define
 $$\rho(\gamma)=\int_{-\infty}^\infty  e^{   \gamma \mathcal{B}_s^\gamma  } Z_s ds.$$
Therefore   the law of the (unit volume or area) quantum sphere $\mu_h$ is obtained as the following limit when $\epsilon$ goes to $0$
\begin{equation}\label{defquantum}
\E[ F(\mu_h)] := \underset{\epsilon \to 0}{\lim}  \; \;   \frac{1}{C_\epsilon}\E\Big[    F\big(  (e^{\star})^2 e^{   \gamma \mathcal{B}_s^\gamma  }  N_\gamma(ds d\theta)      \big)  1_{(e^{\star})^2 \rho(\gamma) \in [1,1+\epsilon] }  \Big]
\end{equation}     
where $C_\epsilon= \E[      1_{(e^{\star})^2 \rho(\gamma) \in [1,1+\epsilon] }  ]$ and  $F$ is an arbitrary bounded measurable functional defined on the space of Radon measures equipped with the weak-$\star$ topology. The asymptotics of this quantity is easily identified with a simple change of variables
\begin{align*}
\E[    F(  (e^{\star})^2 e^{   \gamma \mathcal{B}_s^\gamma  }  N_\gamma(ds d\theta)      )  1_{(e^{\star})^2 \rho(\gamma) \in [1,1+\epsilon] }  ] &= \int_0^\infty \E[    F( u  e^{   \gamma \mathcal{B}_s^\gamma  }  N_\gamma(ds d\theta)  ) 1_{u \rho(\gamma) \in [1,1+\epsilon] }  ]  u^{\frac{\delta}{2}-2} du \\
& = \int_0^\infty \E\Big[  F( y \frac{e^{   \gamma \mathcal{B}_s^\gamma  }  N_\gamma(ds d\theta)}{\rho(\gamma)}  )        \rho(\gamma)^{1-\frac{\delta}{2}}  \Big]  1_{y \in [1,1+\epsilon] }  y^{\frac{\delta}{2}-2} dy .
\end{align*}
It is then straightforward to get the equivalent as $\epsilon\to 0$
$$ \E[    F(  (e^{\star})^2 e^{   \gamma \mathcal{B}_s^\gamma  }  N_\gamma(ds d\theta)      )  1_{(e^{\star})^2 \rho(\gamma) \in [1,1+\epsilon] }  ] \approx  \Big( \int_0^\infty   1_{v \in [1,1+\epsilon] }  v^{\frac{\delta}{2}-2} dv \Big)\E\Big[  F(   \frac{e^{   \gamma \mathcal{B}_s^\gamma  }  N_\gamma(ds d\theta)}{\rho(\gamma)}  )        \rho(\gamma)^{\frac{2}{\gamma}  (Q-\gamma)}  \Big]  $$
where  $a \approx b$ means that the ratio $a/b$ converges to $1$ and we have used that $1-\frac{\delta}{2}= \frac{2}{\gamma}  (Q-\gamma)$. Therefore, definition \eqref{defquantum} is equivalent to  the following  
\begin{definition}
The law of the unit area quantum sphere $\mu_h$ (in the space of quantum surfaces with two marked points $-\infty$ and $\infty$\footnote{see \cite{DMS} for the definition of quantum surfaces. Here, we choose a parameterisation of the unit area quantum sphere such that the radial part is maximal at $0$.}) is given by 
\begin{equation}\label{defquantumsphere} 
\E[  F(\mu_h) ]=   \frac{1}{  \bar{R}(\gamma) }\E\Big[    F(  \frac{e^{   \gamma \mathcal{B}_s^\gamma  }  N_\gamma(ds d\theta)}{\rho(\gamma)} )    \rho(\gamma)^{\frac{2}{\gamma}  (Q-\gamma)}      \Big] 
\end{equation}
where $ \bar{R}(\gamma)$ is the unit volume reflection coefficient and $F$  is an arbitrary bounded measurable functional defined on the space of Radon measures equipped with the weak-$\star$ topology.
\end{definition}
In view of expression \eqref{defquantumsphere}, it is quite natural to interpret $\bar{R}(\gamma)$ as the partition function of the unit area quantum sphere. 
Hence our paper suggests (and proves to the extent of Corollary \ref{th2}) that the unit area quantum sphere is a universal feature  of tails of GMC measures.

Also, the above simple computation shows the equivalence between the reflection coefficient (or 2 point correlation function) constructed in \cite{KRV} and the unit area quantum sphere (with 2 marked points).  Let us mention that this equivalence  \cite{DKRV}$\leftrightarrow$\cite{DMS} has been established in \cite{AHS} for quantum surfaces with 3 marked points.

\subsubsection*{Acknowledgements}  We thank the excellent restaurant {\it Le Pelagos} on the island of Porquerolles where this work was initiated.


\section{Proof of Theorem \ref{th1}}\label{proof}
\subsubsection{Conventions and notations} We will denote $|\cdot|$ the norm in $\C$ of the standard Euclidean   metric. For all $r>0$ we will denote by $B(x,r)$ (resp. $C(x,r)$) the Euclidean ball (resp. circle) of center $x$ and radius $r$. We will write $X\amalg Y$ when the random variables $X$ and $Y $ are independent.

\subsection{Localization trick}\label{sub:localization}
The difficulty in computing right tails of random measures usually lies in the fact that the right tail should depend on the macroscopic shape of the random measure. An important observation in our argument is that tails of GMC measures depend only very softly on their global shape. It is based on the Girsanov transform which we apply in a way so as to  localize the computations for the right tail of GMC measures around a fixed point. Girsanov's transform is a standard tool in the study of Gaussian processes (or even GMC measures) but, to the best of our knowledge, it has never been used  to localize large mass effects of GMC measures. Anyway, this simple observation reduces drastically the difficulty of computing the tails. 
\begin{lemma}\label{localization}{ (Localization trick)}
Let $\mathcal{O}$ be an open set. We have 
\begin{align}\label{magic}
\P(M_\gamma(\mathcal{O})>t)&=  \int_\mathcal{O}\E[\frac{1}{M_\gamma(v,\mathcal{O}) } \mathbf{1}_{\{M_\gamma(v,\mathcal{O})>t\}}] d^2v
\end{align}
where we have set 
$$M_\gamma(v,\mathcal{O}):=\int_{\mathcal{O}}\frac{1}{|z-v|^{\gamma^2}}M_\gamma(d^2z).$$
\end{lemma}

\noindent {\it Proof.}
First write
\begin{align}\label{girs}
\P(M_\gamma(\mathcal{O})>t)&=  \E[\frac{M_\gamma(\mathcal{O})}{M_\gamma(\mathcal{O})}\mathbf{1}_{\{M_\gamma(\mathcal{O})>t\}}]=\int_{\mathcal{O}}  \E[e^{\gamma X(v)-\frac{\gamma^2}{2}\E[X(v)^2]}\frac{1}{M_\gamma(\mathcal{O})}\mathbf{1}_{\{M_\gamma(\mathcal{O})>t\}}]\,d^2v.
\end{align}
Girsanov's transform then asserts that weighting the probability law by $e^{\gamma X(v)-\frac{\gamma^2}{2}\E[X(v)^2]}$ amounts to shifting the law of $X$ by $\gamma \ln\frac{1}{|\cdot-v|}$ and our claim follows.
Of course, as the field $X$ is distribution-valued, it cannot be evaluated pointwise: the last equality in \eqref{girs} thus requires a cut-off regularization of $X(v)$ and then a harmless passage to the limit in the regularization parameter to be rigorously established (use  the fact that $\sup_{\epsilon>0}\E[\frac{1}{M_{\gamma,\epsilon}(\mathcal{O})}]<+\infty$ where $\epsilon$ stands for the regularization parameter).   \qed
 
\medskip

The important point is that the asymptotic behavior  of the quantity $\E[\frac{1}{M_\gamma(v,\mathcal{O}) } \mathbf{1}_{\{M_\gamma(v,\mathcal{O})>t\}}]$ is completely dominated  by the behavior close to $v$ because of the singularity $\frac{1}{|z-v|^{\gamma^2}}$. So we will establish the following estimate, which will be enough to complete our argument.

\begin{lemma}\label{main}
For any $\delta\in (0,\frac{1+p_0-\frac{4}{\gamma^2}}{2+p_0})$ and for all open sets  $\mathcal{O}\subset D$ with $C^1$ boundary, there exists some function $\epsilon:\R_+\times \mathcal{O}\to\R$ such that 
$$\forall v \in \mathcal{O},\quad  \E[\frac{1}{M_\gamma(v,\mathcal{O}) } \mathbf{1}_{\{M_\gamma(v,\mathcal{O})>t\}}]=   \frac{\frac{2}{\gamma}(Q-\gamma)}{\frac{2}{\gamma}(Q-\gamma)+1}   \frac{\bar{R}(\gamma)}{t^{\frac{4}{\gamma^2}}}+ \frac{1}{t^{\frac{4}{\gamma^2}+\delta}}\epsilon(t,v)$$
and satisfying $\forall  v$ $\lim_{t\to\infty} |\epsilon(t,v)|=0$ and $  \sup_{v\in \mathcal{O},t\geq 1}{\rm dist}(v,\mathcal{O}^c)^\alpha|\epsilon(t,v)|<+\infty$ for some $\alpha<1$.
\end{lemma}

Now, relation \eqref{girs} combined with Lemma \ref{main} yields that 
\begin{equation*}
\P(M_\gamma(\mathcal{O})>t)=  \frac{\frac{2}{\gamma}(Q-\gamma)}{\frac{2}{\gamma}(Q-\gamma)+1}   \frac{\bar{R}(\gamma)}{t^{\frac{4}{\gamma^2}}}+ \frac{1}{t^{\frac{4}{\gamma^2}+\delta}} \int_{\mathcal{O}} \epsilon(t,v) \,d^2v. 
\end{equation*}
Notice that since $\mathcal{O}$ has a $C^1$ boundary and $ \sup_{v\in \mathcal{O},t\geq 1}{\rm dist}(v,\mathcal{O}^c)^\alpha|\epsilon(t,v)|<+\infty$, the term  $\int_{\mathcal{O}} \epsilon(t,v) \,d^2v$ is indeed well defined. In fact, the term $\int_{\mathcal{O}} \epsilon(t,v) \,d^2v$ converges to $0$ as $t$ goes to infinity by using the dominated convergence theorem and the fact that for all $v$ one has $\lim_{t\to\infty} |\epsilon(t,v)|=0$. This yields Theorem \ref{th1}. The remaining part of this section is thus devoted to the proof of this lemma.

\begin{remark}\label{C1assump}
Let us stress here that the assumption that $\mathcal{O}$ has a $C^1$ boundary ensures that   $  \int_{\mathcal{O}}  \frac{1}{{\rm dist}(v,\mathcal{O}^c)^\alpha}\,d^2v < \infty$ for all $\alpha<1$.  Therefore the term $\int_{\mathcal{O}} \epsilon(t,v) \,d^2v$ is well defined and converges to $0$. This is the only place where the $C^1$ assumption is used and perhaps one could relax a bit the $C^1$ assumption by a finer analysis of the $\epsilon(t,v)$ term. 
\end{remark}
 
\subsection{GFF with vanishing mean over a circle}
Given $v\in \C$, we consider  the GFF $X^{v,r}$ with vanishing mean over the circle $C(v,r)$, namely  a centered Gaussian random distribution with covariance  for $x,y\in D$
\begin{equation}\label{hatGformula}
\E [ X^{v,r}(x)X^{v,r}(y)] 
=\ln\frac{1}{|x-y|}+\ln\big(\frac{|x-v|}{r}\big)_++\ln\big(\frac{|y-v|}{r}\big)_+ +\ln r
\end{equation}
where we use the notation $ |z|_+=|z|$ if $|z|\geq 1$ and  $ |z|_+=1$ if $|z|\leq 1$. 
 
Let us denote by $(X^{v,r}_u )_{u\geq 0}$ the circle average of this field 
    $$X^{v,r}_u=\frac{1}{2\pi i}\oint_{|w-v|={re^{-u}}} X^{v,r}(w)\frac{dw}{w-v}.$$
  A simple computation shows that $(X^{v,r}_u)_{u \geq 0}$ is a standard Brownian motion starting from $0$ independent of the sigma algebra $\sigma\{X^{v,r}(z);z\in B(v,r)^c\}$.

%
%
%
  
  \subsection{Polar decomposition and the reflection coefficient} \label{Reflcoefficient}  
  
  The asymptotic expansion in Lemma \ref{main}   will be determined by the behavior of the integral $M_\gamma(v,\mathcal{O})$ around the singularity at $v$. As already noticed in \cite{KRV}, an important ingredient in the analysis is the reflection coefficient. We consider the polar decomposition of the chaos measure.  We have the following equality in the sense of distributions
\begin{equation*}
X^{v,r}(v+re^{-s} e^{i \theta})=B_s+Y(s,\theta)
\end{equation*}
where $B_s$ is Brownian Motion starting form the origin at $s=0$ and $Y(s,\theta)$ the independent centered Gaussian field with  covariance 
\begin{equation}\label{covlateral}
\E[  Y(s,\theta) Y(s',\theta') ]  = \ln \frac{e^{-s}\vee e^{-s'}}{|e^{-s}e^{i \theta} - e^{-s'} e^{i \theta'} |}.
\end{equation}
Then we get by the change of variables $z=v+re^{-s}e^{i\theta}$\footnote{Observe that the  measures involved in this relation are defined through a limiting procedure. Thus one needs to apply the change of variables in the regularized measures and then pass to the limit to obtain the relation.}
\begin{equation}\label{ialphadef}
\int_{B(v,r)}\frac{e^{\gamma X^{v,r}(z)-\frac{\gamma^2}{2} \E[   X^{v,r}(z)^2 ] }}{|z-v|^{\gamma^2}}d^2z\stackrel{law}=r^{2-\gamma^2}\int_0^\infty  e^{\gamma (B_s-(Q-\gamma)s)}Z_s ds 
\end{equation}

The following decomposition lemma due to Williams (see \cite{Williams})  will be useful in the study of $M_\gamma(v,B(v,r))$:

\begin{lemma}\label{lemmaWilliams}
Let $(B_s-\nu s)_{s \geq 0}$ be a Brownian motion with negative drift, i.e. $\nu >0$ and let $M= \sup_s (B_s-\nu s)$.  Then conditionally on $M$ the law of the path   $(B_s-\nu s)_{s \geq 0}$ is given by the joining of two  independent paths:

\begin{itemize}
\item
A Brownian motion $((B_s^1+\nu s))_{s \leq \tau_M}$ with positive drift $\nu >0$ run until its hitting time $\tau_M$ of $M$. 
\item
$(M+{B}^2_s- \nu s)_{s\geq 0}$ where ${B}^2_s- \nu s$  is a  Brownian motion with negative drift conditioned to stay negative.
\end{itemize}
By the joining of two paths $(X_s)_{s \leq T}$ (with $T>0$) and $(\bar{X}_s)_{s \geq 0}$, we mean the path $(X_s 1_{s \leq T}+\bar{X}_{s-T} 1_{s > T} )_{s \geq 0}$. Moreover, one has the following time reversal property for all $C>0$ (where $\tau_C$ denotes the hitting time of $C$) 
\begin{equation*}
(B_{\tau_C-s}^1+\nu (\tau_C-s)-C)_{s  \leq \tau_C}\stackrel{law}=  (\tilde{B}_s- \nu s)_{s \leq L_{-C}}
\end{equation*}
where $(\tilde{B}_s- \nu s)_{s \geq 0}$ is a Brownian motion with drift $-\nu$ conditioned to stay negative and $ L_{-C}$ is the last time  $(\tilde{B}_s- \nu s)$ hits $-C$.  

\end{lemma}

\begin{remark}\label{fundamentalremark}
As a consequence of the above lemma, one can also deduce that the process $(\tilde{B}_{L_{-C}+s}- \nu (L_{-C}+s)+C)_{s \geq 0}$ is equal in distribution to $(\tilde{B}_s- \nu s)_{s \geq 0}$. 
\end{remark}
 We may apply Lemma \ref{lemmaWilliams} to \eqref{ialphadef}. Let $M=\sup_{s\geq 0}(B_s-(Q-\gamma)s)$ and $L_{-M}$ be the last time $(B^\gamma_s)_{s \geq 0}$ hits $-M$ (recall that $B^\gamma_s$ was defined in \eqref{BMneg}). Then
\begin{align}
&\int_0^\infty e^{\gamma (B_s-(Q-\gamma)s)}Z_s ds\stackrel{law}=e^{\gamma M}\int_{-L_{-M}}^\infty e^{\gamma\caB_s^{\gamma}}Z_{s+L_{-M}}ds\stackrel{law}=e^{\gamma M}\int_{-L_{-M}}^\infty e^{\gamma\caB_s^{\gamma}}Z_{s} ds\label{will11}
\end{align}
where we used stationarity of the process $Z_s$ (and independence of $Z_s$ and $B_s$). The distribution of $M$ is well known (see section 3.5.C in the textbook \cite{KaraSh} for instance):
\begin{equation}\label{tailofmax}
 \P (e^{\gamma M}>t) = \frac{1}{t^{\frac{2(Q-\gamma)}{\gamma}}} ,  \quad  t\geq 1.
\end{equation}
Now, notice that $\bar{R}(\gamma)$ appears as the coefficient of the tail of the random variable $e^{\gamma M}\int_{-\infty}^\infty e^{\gamma\caB_s^{\gamma}}Z_{s} ds$
$$\forall \eta>0,\quad \P\big(e^{\gamma M}\int_{-\infty}^\infty e^{\gamma\caB_s^{\gamma}}Z_{s} ds>t\big) = \bar{R}(\gamma)t^{-\frac{2}{\gamma}(Q-\gamma)}+o(t^{-\frac{2}{\gamma}(Q-\gamma)-(1-\eta)})$$
Indeed, setting $\mathcal{Z}:=\int_{-\infty}^\infty e^{\gamma\caB_s^{\gamma}}Z_{s} ds$, we have because of \eqref{tailofmax} 
  
\begin{align*}\P\big(e^{\gamma M}\mathcal{Z}>t\big)=&\P\big(e^{\gamma M}\mathcal{Z}>t,\mathcal{Z}\leq t\big)+\P\big(e^{\gamma M}\mathcal{Z}>t,\mathcal{Z}> t\big)\\
=&\E[\mathcal{Z}^{\frac{2}{\gamma}(Q-\gamma)}\mathbf{1}_{\mathcal{Z}\leq t}]t^{-\frac{2}{\gamma}(Q-\gamma)}+\P\big(e^{\gamma M}\mathcal{Z}>t,\mathcal{Z}> t\big)\\
 =&\bar{R}(\gamma)t^{-\frac{2}{\gamma}(Q-\gamma)}+ \E[\mathcal{Z}^{\frac{2}{\gamma}(Q-\gamma)}\mathbf{1}_{\mathcal{Z}> t}]t^{-\frac{2}{\gamma}(Q-\gamma)}+\P\big(e^{\gamma M}\mathcal{Z}>t,\mathcal{Z}> t\big)
\end{align*} 
Now, using  \eqref{momentcondR} and Markov's inequality, the two correction terms are easily seen to be $o(t^{-\frac{2}{\gamma}(Q-\gamma)-(1-\eta)})$.
%
%

\subsection{Decomposition of the GMC around $v$}
We fix  $v\in \mathcal{O}$ and set $ r={\rm dist}(v,\partial \mathcal{O})$.  Our purpose now is to replace the field $X$ in the definition of $M_\gamma(v,B(v,r))$ by the field $X^{v,r}$ in order to use the polar decomposition of the GMC measure in $B(v,r)$ detailed in the previous subsection. It suffices to subtract the mean value of $X$ along the circle $C(v,r)$. Therefore we introduce
\begin{equation*}
N_{v,r}=\frac{1}{2\pi i} \oint_{|w-v|=r} X  (w) \frac{dw}{w-v}.  
\end{equation*}
This is a centered Gaussian variable with variance
\begin{equation*}
\E[N_{v,r}^2]=- \ln r.
\end{equation*}
We introduce the field $\tilde{X}= X -N_{v,r}$, which has the same law as $X^{v,r}$. One can notice that 
\begin{equation*}
\E[  \tilde{X}(x) N_{v,r}  ]= 0
\end{equation*}
for all $|x-v| \leq r$ hence $\tilde{X}$ is independent from $N_{v,r}$ in this region. Therefore, using  the decomposition in distribution obtained in subsection \ref{Reflcoefficient}  
\begin{equation}\label{decomp}
M_\gamma(v,B(v,r))= r^{2-\gamma^2} e^{\gamma M} E_{r,v}  I_\gamma(M)
\end{equation}
where we have set
\begin{align*}
E_{r,v} &:=e^{\gamma N_{v,r}  -\frac{\gamma^2}{2} \E[ N_{v,r}^2   ] }, & I_\gamma(M)&:= \int_{-L_{-M}}^\infty e^{\gamma\caB_s^{\gamma}}Z_{s} ds
\end{align*}
 with the convention that $I_\gamma(\infty)= \int_{-\infty}^\infty e^{\gamma\caB_s^{\gamma}}Z_{s} ds$. We will constantly use in the sequel the fact that $x \mapsto I_\gamma(x)$ is an increasing function.

\subsection{Getting rid of the non-singularity}\label{rid}
Now we argue that the behaviour of $\E[\frac{1}{M_\gamma(v,\mathcal{O}) } \mathbf{1}_{\{M_\gamma(v,\mathcal{O})>t\}}]$ is completely determined by that of $\E[\frac{1}{M_\gamma(v,B(v,r))} \mathbf{1}_{\{M_\gamma(v,B(v,r))>t\}}]$. We will show this claim by giving an upper/lower bound of $\E[\frac{1}{M_\gamma(v,\mathcal{O}) } \mathbf{1}_{\{M_\gamma(v,\mathcal{O})>t\}}]$ up to $o(t^{-\frac{4}{\gamma^2}-\delta})$ for any $\delta\in (0,\frac{1+p_0-\frac{4}{\gamma^2}}{2+p_0})$. We introduce the notation $A:=M_\gamma(v,B(v,r)^c\cap\mathcal{O})$. A key lemma in the argument is the following.
\begin{lemma}\label{lemA}
For all $p \in (0,\frac{4}{\gamma^2})$ there is  some constant $C_p>0$,  
$$ \forall v\in \mathcal{O},\quad \E[A^p]\leq C_p \times \left\{\begin{array}{ll}1 & \text { if }p \in (0,\frac{4}{\gamma^2}-1) \\
(\ln \frac{1}{r})^{p\vee 1} &  \text { if }p=\frac{4}{\gamma^2}-1\\
 r^{\psi(p )}&  \text { if }\frac{4}{\gamma^2}-1<p<\frac{4}{\gamma^2}\\
\end{array}\right.$$  
with $\psi(p)=(2-\frac{\gamma^2}{2})p-\frac{\gamma^2}{2}p^2$ and $r={\rm dist}(v,\mathcal{O}^c)$. 
\end{lemma}
\begin{remark}\label{spectrum}
Notice that $\psi(p)>0$ for $p\in (0,4/\gamma^2-1)$. The value of $p_0$ in \eqref{defp0}  is determined in such a way that  $\psi(p_0)=-1$. This ensures that $\psi(p)\in (-1,0)$ for $p\in (4/\gamma^2-1,p_0)$. Our proof will show that the exponent $\alpha$ appearing in Lemma \ref{main} will actually be given by $\psi(p)$ for some $p\in (4/\gamma^2-1,p_0)$.
\end{remark}
\noindent {\it Proof.}  
If $p>0$ observe that $\psi(p)> 0 $ if and only if $p<\frac{2}{\gamma}(Q-\gamma)=\frac{4}{\gamma^2}-1$. In this case, the result is a simple consequence of Lemma A.1 in \cite{DKRV}.\\
Now we study the case $p\geq \frac{4}{\gamma^2}-1$. Set $B_n :=\{x\in D;2^{-n-1}\leq |x-v|\leq 2^{-n}\}$ and $A_n=M_\gamma(0,B_n)^p$. Assume $p>1$ (if $p\leq 1$ use sub-additivity of the mapping $x\mapsto x^p$ in the next argument). Using in turn Minkowski's inequality and  invariance under translations we have
$$ \E[A^p]^{1/p}\leq \sum_{n=0}^{\frac{-\ln r}{\ln 2}}\E[A_n^p]^{1/p}+\E[M_\gamma(\mathcal{O})^p]^{1/p}.$$
It is standard fact in GMC theory that $A_n$ scales as $\E[A_n^p]=2^{-n\psi(p)}\E[A_1^p]$ (see the review \cite{review} for instance). We deduce
that
$$\E[A^p]^{1/p}\leq C_p \sum_{n=0}^{\frac{-\ln r}{\ln 2}}2^{-n\psi(p)/p}.$$
Since $\psi(p)\leq 0$, we get our claim.

\qed

\medskip
Now we first give the upper bound. Fix $\delta\in (0,\frac{p_0+1-\tfrac{4}{\gamma^2}}{2+p_0})$. We can find $\eta>\delta$ such that  
\begin{equation}\label{delta}
(1-\eta)(1+p_0)> \frac{4}{\gamma^2}+\delta, \quad 1+p_0> \frac{4}{\gamma^2}+\eta.
\end{equation}
Hence we can choose $\frac{4}{\gamma^2}<p<p_0+1$ such that $(1-\eta)p>\frac{4}{\gamma^2}+\delta$.
\begin{align}
\E[\frac{1}{M_\gamma(v,\mathcal{O}) } \mathbf{1}_{\{M_\gamma(v,\mathcal{O})>t\}}]\leq & \E[\frac{1}{M_\gamma(v,B(v,r))+A } \mathbf{1}_{\{M_\gamma(v,B(v,r))>t-t^{1-\eta}\}}]+\E[\frac{1}{A } \mathbf{1}_{\{A>t^{1-\eta}\}}] \nonumber\\
\leq & \E[\frac{1}{M_\gamma(v,B(v,r))} \mathbf{1}_{\{M_\gamma(v,B(v,r))>t-t^{1-\eta}\}}]+\E[A^{p-1}]t^{-(1-\eta)p}.  \label{upperbound}
\end{align}
Notice that the expectation $\E[A^{p-1}]$ can be analyzed with Lemma \ref{lemA} and gives $\E[A^{p-1}]\leq C_pr^{\psi(p-1)}$ with $\psi(p-1)\in (-1,0)$.
 
For the lower bound we first restrict to the set $\{A<t^{1-\eta}\}$  to get  
\begin{align}
 \E[ &\frac{1}{M_\gamma(v,B(v,r))+A } \mathbf{1}_{\{M_\gamma(v,B(v,r))+A>t\}}] \nonumber \\
\geq &\E[\frac{1}{M_\gamma(v,B(v,r)) +t^{1-\eta}} \mathbf{1}_{\{M_\gamma(v,B(v,r))>t,A<t^{1-\eta}\}}] \nonumber  \\
= &\E[\frac{1}{M_\gamma(v,B(v,r)) +t^{1-\eta}} \mathbf{1}_{\{M_\gamma(v,B(v,r))>t\}}]-\E[\frac{1}{M_\gamma(v,B(v,r)) +t^{1-\eta}} \mathbf{1}_{\{M_\gamma(v,B(v,r))>t,A>t^{1-\eta}\}}] \nonumber  \\
\geq &\E[\frac{1}{M_\gamma(v,B(v,r)) } \mathbf{1}_{\{M_\gamma(v,B(v,r))>t\}}]-t^{1-\eta}\E[\frac{1}{M_\gamma(v,B(v,r))^2 } \mathbf{1}_{\{M_\gamma(v,B(v,r))>t\}}] \nonumber  \\
&-\E[\frac{1}{M_\gamma(v,B(v,r)) +t^{1-\eta}} \mathbf{1}_{\{M_\gamma(v,B(v,r))>t,A>t^{1-\eta}\}}]   \nonumber  \\
\geq &  (1-\frac{1}{t^\eta})\E[\frac{1}{M_\gamma(v,B(v,r)) } \mathbf{1}_{\{M_\gamma(v,B(v,r))>t\}}]   -\E[\frac{1}{M_\gamma(v,B(v,r)) +t^{1-\eta}} \mathbf{1}_{\{M_\gamma(v,B(v,r))>t,A>t^{1-\eta}\}}]   \label{lowerbound}
\end{align}
where we have used the inequality $(1+u)^{-1}\geq 1-u$. The last of the three terms $$\E[\frac{1}{M_\gamma(v,B(v,r)) +t^{1-\eta}} \mathbf{1}_{\{M_\gamma(v,B(v,r))>t,A>t^{1-\eta}\}}]$$ is less than $t^{-(1-\eta)}\P(A>t^{1-\eta})$ and by the Markov inequality, 
\begin{equation*}
t^{-(1-\eta)}\P(A>t^{1-\eta}) \leq t^{-(1-\eta)p} \E[A^{p-1}] \leq C  \leq t^{-(1-\eta)p} r^{\psi(p-1)}
\end{equation*}
with $\psi(p-1)\in (-1,0)$. The second will be treated in the same time as the first (next subsection) and will be proved to be $O(t^{-\frac{4}{\gamma^2}-\eta})$, hence $o(t^{-\frac{4}{\gamma^2}-\delta})$.

\subsection{Behaviour near the singularity}\label{nearsing}

In view of the bounds  \eqref{upperbound} and \eqref{lowerbound}, it remains to study the term 
\begin{equation*}
\E[\frac{1}{M_\gamma(v,B(v,r)) } \mathbf{1}_{\{M_\gamma(v,B(v,r))>t\}}].
\end{equation*}

We will proceed by  using  decomposition \eqref{decomp} and  establishing lower/upper bounds. 

\subsubsection*{Lower bound}

In the study of the lower bound, let us fix $\eta$ such that $\eta>\delta$ and $\eta$ satisfies \eqref{delta} in such a way that we can find $1<p<1+p_0$ with $p(1-\eta)>\frac{4}{\gamma^2}+\delta$. We will introduce the event $\{\gamma M> \eta\ln t \}$. On this event, one has $I_\gamma (M) \geq I_\gamma( \frac{ \eta}{\gamma}\ln t )$ and therefore 
\begin{align*}
& \E[\frac{1}{M_\gamma(v,B(v,r))} \mathbf{1}_{\{M_\gamma(v,B(v,r))>t\}}]  \\
&= \E[\frac{1}{r^{2-\gamma^2} E_{r,v}  e^{\gamma M}I_\gamma(M)} \mathbf{1}_{\{r^{2-\gamma^2} E_{r,v} e^{\gamma M}I_\gamma( M)>t\}} ] \\
&\geq \E[\frac{1}{r^{2-\gamma^2} E_{r,v}  e^{\gamma M}I_\gamma(\infty)} \mathbf{1}_{\{r^{2-\gamma^2} E_{r,v} e^{\gamma M}I_\gamma( \frac{\eta}{\gamma}\ln t )>t\}} \mathbf{1}_{\{\gamma M> \eta\ln t \}}] \\
&= \E[\frac{1}{r^{2-\gamma^2} E_{r,v}  e^{\gamma M}I_\gamma(\infty)} \mathbf{1}_{\{r^{2-\gamma^2} E_{r,v} e^{\gamma M}I_\gamma( \frac{ \eta}{\gamma}\ln t )>t\}} ] \\
&\,\,\,  -\E[\frac{1}{r^{2-\gamma^2} E_{r,v}  e^{\gamma M}I_\gamma(\infty)} \mathbf{1}_{\{r^{2-\gamma^2} E_{r,v} e^{\gamma M}I_\gamma( \frac{ \eta}{\gamma}\ln t )>t\}} \mathbf{1}_{\{\gamma M< \eta\ln t \}}] 
\end{align*}
The correction term (second expectation in final line above) is $r^{-\alpha}o(t^{-\frac{4}{\gamma^2}-\delta})$ with $\alpha<1$. Indeed on the event $\{\gamma M<\eta\ln t \}$, the event  $\{r^{2-\gamma^2} E_{r,v} e^{\gamma M}I_\gamma( \frac{ \eta}{\gamma}\ln t )>t\}$ is contained in the event $\{r^{2-\gamma^2} E_{r,v} I_\gamma( \infty )>t^{1-\eta}\}$. Hence the correction term is less than $\E[(r^{2-\gamma^2} E_{r,v} )^{p-1}]t^{-(1-\eta)p}\E[I_\gamma(\infty)^{p-1}]=r^{\psi(p-1)}t^{-(1-\eta)p}\E[I_\gamma(\infty)^{p-1}]$, where the function $\psi$ has been defined in Lemma \ref{lemA}.

Then we average with respect to $M$ to get (recall that the distribution of $M$ is given by \eqref{tailofmax}) 
\begin{align*}
& \E[\frac{1}{r^{2-\gamma^2} E_{r,v}  e^{\gamma M}I_\gamma(\infty)} \mathbf{1}_{\{r^{2-\gamma^2} E_{r,v} e^{\gamma M}I_\gamma( \frac{ \eta}{\gamma}\ln t )>t\}} ] \\
& \geq \E[\frac{1}{r^{2-\gamma^2} E_{r,v}  e^{\gamma M}I_\gamma(\infty)} \mathbf{1}_{\{r^{2-\gamma^2} E_{r,v} e^{\gamma M}I_\gamma( \frac{ \eta}{\gamma}\ln t )>t\}} \mathbf{1}_{\{r^{2-\gamma^2} E_{r,v} I_\gamma( \infty)\leq t\}} ] \\
&= \frac{2(Q-\gamma)}{2(Q-\gamma)+\gamma}
 \E[ \big (   r^{2-\gamma^2}E_{r,v} \big )^{\frac{2}{\gamma}(Q-\gamma)}  \frac{I_\gamma( \frac{ \eta }{\gamma}\ln t )^{\frac{4}{\gamma^2}}}{I_\gamma(\infty)} \mathbf{1}_{\{r^{2-\gamma^2} E_{r,v} I_\gamma( \infty)\leq t\}}]t^{-\frac{4}{\gamma^2}}\\
&= \frac{2(Q-\gamma)}{2(Q-\gamma)+\gamma}
 \E\big[\frac{I_\gamma( \frac{ \eta }{\gamma}\ln t )^{\frac{4}{\gamma^2}}}{I_\gamma(\infty)}\big]t^{-\frac{4}{\gamma^2}} \\
&- \frac{2(Q-\gamma)}{2(Q-\gamma)+\gamma}
 \E[ \big (   r^{2-\gamma^2}E_{r,v} \big )^{\frac{2}{\gamma}(Q-\gamma)}  \frac{I_\gamma( \frac{ \eta }{\gamma}\ln t )^{\frac{4}{\gamma^2}}}{I_\gamma(\infty)} \mathbf{1}_{\{r^{2-\gamma^2} E_{r,v} I_\gamma( \infty)> t\}}]t^{-\frac{4}{\gamma^2}}.
\end{align*}
Above we have used the fact $ \E[ \big (   r^{2-\gamma^2}E_{r,v} \big )^{\frac{2}{\gamma}(Q-\gamma)} ]=1$. One has 
\begin{align*}
 \E[ \big (   r^{2-\gamma^2}E_{r,v} \big )^{\frac{2}{\gamma}(Q-\gamma)}  \frac{I_\gamma( \frac{ \eta }{\gamma}\ln t )^{\frac{4}{\gamma^2}}}{I_\gamma(\infty)} \mathbf{1}_{\{r^{2-\gamma^2} E_{r,v} I_\gamma( \infty)> t\}}]t^{-\frac{4}{\gamma^2}}
& \leq \E[ \big (   r^{2-\gamma^2}E_{r,v} \big )^{\frac{4}{\gamma^2}-1+\eta}  I_\gamma( \infty)^{\frac{4}{\gamma^2}-1+\eta}]t^{-\frac{4}{\gamma^2}-\eta}  \\
&=  r^{\psi (\frac{4}{\gamma^2}-1+\eta)} \E[   I_\gamma( \infty)^{\frac{4}{\gamma^2}-1+\eta} ] t^{-\frac{4}{\gamma^2}-\eta}  \\
\end{align*}
where $\frac{4}{\gamma^2}-1+\eta<p_0$ by condition \eqref{delta}. Hence it remains to show that $  \E\big[\frac{I_\gamma( \frac{\eta }{\gamma}\ln t )^{\frac{4}{\gamma^2}}}{I_\gamma(\infty)}\big]=\bar{R}(\gamma)(1+o(t^{-\eta }))$ to get the desired lower bound. By Remark \ref{fundamentalremark}, the process $\hat{B}_s^{\gamma}$ defined for $s \leq 0$ by the relation $\hat{B}^{\gamma}_s= \mathcal{B}^{\gamma}_{s-L_{-\frac{\eta }{\gamma} \ln t} }+ \frac{\eta }{\gamma} \ln t$ is independent from everything and distributed like $(\mathcal{B}^{\gamma}_s)_{s \leq 0}$. We can then write
\begin{equation*}
 \int_{-\infty}^\infty  e^{\gamma \mathcal{B}^{\gamma}_s} Z_s ds  =I_\gamma( \frac{\eta }{\gamma}\ln t ) +{t^{-\eta }}B
\end{equation*}
with
\begin{equation*}
B= \int_{-\infty}^0 e^{\gamma \hat{B}^{\gamma}_s}Z_{s- L_{-\frac{\eta}{\gamma} \ln t}} ds.
\end{equation*}
We set $m=\frac{4}{\gamma^2}$. Then we have (use the triangle inequality to get the third line below)
\begin{align*}
 \bar{R}&(\gamma)- \E[I_\gamma( \frac{ \eta }{\gamma}\ln t )^{m}/I_\gamma(\infty)]\\
 =&  \E[\Big(\frac{I_\gamma( \frac{\eta }{\gamma}\ln t )+t^{-\eta  }B}{(I_\gamma( \frac{ \eta}{\gamma}\ln t )+t^{-\eta }B)^{1/m}}\Big)^m]-\E[\Big(\frac{I_\gamma( \frac{ \eta }{\gamma}\ln t ) }{(I_\gamma( \frac{\eta }{\gamma}\ln t )+t^{-\eta }B)^{1/m}}\Big)^m]\\
\leq & \Big(\E[\Big(\frac{I_\gamma( \frac{ \eta}{\gamma}\ln t )}{(I_\gamma( \frac{\eta}{\gamma}\ln t )+t^{-\eta }B)^{1/m}}\Big)^m]^{1/m}+\E[\Big(\frac{ t^{-\eta  }B}{(I_\gamma( \frac{\eta}{\gamma}\ln t )+t^{-\eta}B)^{1/m}}\Big)^m]^{1/m}\Big)^m-\E[\Big(\frac{I_\gamma( \frac{ \eta }{\gamma}\ln t ) }{(I_\gamma( \frac{ \eta }{\gamma}\ln t )+t^{-\eta }B)^{1/m}}\Big)^m].
\end{align*}  
This expression is of the form $(\alpha_t+\beta_t)^m-\alpha_t^m$ with $\frac{\beta_t}{\alpha_t}\leq c$ for some $c>0$ and all $t\geq 0$. Let us consider another constant $C$ such that $(1+x)^m-1\leq Cx$ for $0\leq x\leq c$. We deduce that  $(\alpha_t+\beta_t)^m-\alpha_t^m\leq C \alpha_t^{m-1} \beta_t$. Plugging this estimate in the above expression yields

\begin{align*}
 \bar{R}&(\gamma)- \E[I_\gamma( \frac{ \eta }{\gamma}\ln t )^{m}/I_\gamma(\infty)]\\
\leq &C t^{-\eta } \E[\Big(\frac{I_\gamma( \frac{ \eta}{\gamma}\ln t )}{(I_\gamma( \frac{\eta}{\gamma}\ln t )+t^{-\eta }B)^{1/m}}\Big)^m]^{1-1/m}\E[\Big(\frac{ B}{(I_\gamma( \frac{ \eta }{\gamma}\ln t )+t^{-\eta }B)^{1/m}}\Big)^m]^{1/m}.
\end{align*}  
In the above expression the first expectation is less  than $\bar{R}(\gamma)^{1-\gamma^2/4}$. We aim to bound the second one. Let us fix $\eta'\in (0,1)$. Then 
$$\E[\Big(\frac{ B}{(I_\gamma( \frac{ \eta }{\gamma}\ln t )+t^{-\eta  }B)^{1/m}}\Big)^m]^{1/m}\leq \E[ \frac{ B^m}{I_\gamma(0 )^{1-\eta'}(t^{-\eta  }B)^{\eta'}} ]^{1/m}\leq t^{\frac{\eta' \eta }{m}}\E[B^{m-\eta'}   I_\gamma( 0)^{-1+\eta'}]^{1/m}.$$
Notice that  the  last expectation is  finite (use H\"{o}lder inequality and the fact that $B$ has finite moments of order $q<m$ and $ I_\gamma( 0)$ has finite negative moments of all order).  This proves our claim for the lower bound provided that we choose $\eta'$ small enough so as to make $\eta(1-\frac{\eta'}{m})>\delta$.
  
\subsubsection*{Upper bound}
  
For the upper bound, we use again the decomposition \eqref{decomp} 
\begin{align*}
 \E[\frac{1}{M_\gamma(v,B(v,r))} \mathbf{1}_{\{M_\gamma(v,B(v,r))>t\}}]  
&=\E[\frac{1}{r^{2-\gamma^2} E_{r,v}  e^{\gamma M}I_\gamma(M)} \mathbf{1}_{\{r^{2-\gamma^2} E_{r,v} e^{\gamma M}I_\gamma( M )>t\}} ]  
\end{align*}
We want to replace the term $I_\gamma(M)$ in the fraction by $I_\gamma(\infty)$. Hence we write 
\begin{align*}
 \E[\frac{1}{M_\gamma(v,B(v,r))} \mathbf{1}_{\{M_\gamma(v,B(v,r))>t\}}]  
&=\E[\frac{1}{r^{2-\gamma^2} E_{r,v}  e^{\gamma M}I_\gamma(\infty)} \mathbf{1}_{\{r^{2-\gamma^2} E_{r,v} e^{\gamma M}I_\gamma( M )>t\}} ] +C(t)
\end{align*}
where $C(t)$ stands for the  cost for this replacement
\begin{align}\label{cost}
C(t):=\E[&\frac{1}{r^{2-\gamma^2}  E_{r,v}  e^{\gamma M}}\big(\frac{1}{I_\gamma(M)}-\frac{1}{I_\gamma(\infty)}\big) \mathbf{1}_{\{r^{2-\gamma^2} E_{r,v} e^{\gamma M}I_\gamma( M )>t\}} ]. 
\end{align}
Now we establish that this cost satisfies  $C(t)=r^{-\alpha}o(t^{-\frac{4}{\gamma^2}-\delta})$ with $\alpha<1$.  

By Remark \ref{fundamentalremark}, conditionally on $M$, the process $\hat{B}_s^{\gamma}$ defined for $s \leq 0$ by the relation $\hat{B}^{\gamma}_s= \mathcal{B}^{\gamma}_{s-L_{-M} }+M$ is independent from everything and distributed like $(\mathcal{B}^{\gamma}_s)_{s \leq 0}$. We can then write
\begin{equation*}
 \int_{-\infty}^\infty  e^{\gamma \mathcal{B}^{\gamma}_s} Z_s ds  =I_\gamma( M) +e^{-\gamma M}B
\end{equation*}
with
\begin{equation*}
B= \int_{-\infty}^0 e^{\gamma \hat{B}^{\gamma}_s}Z_{s- L_{-M}} ds.
\end{equation*}
 Now we observe that, under the condition $\delta\in (0,\frac{p_0+1-\tfrac{4}{\gamma^2}}{2+p_0})$ one can find $\eta\in(0,1)$ and $p\in (0,1+p_0)$ such that
\begin{equation}\label{condition}
p(1-\eta)+\frac{4}{\gamma^2}\eta>\delta+\frac{4}{\gamma^2}\quad \text{ and }\quad \eta(1+\frac{4}{\gamma^2})>\delta+\frac{4}{\gamma^2}.
\end{equation}
Indeed this condition is equivalent to the set of conditions $(1+p_0)(1-\eta)+\frac{4}{\gamma^2}\eta>\delta+\frac{4}{\gamma^2}$ and $\eta(1+\frac{4}{\gamma^2})>\delta+\frac{4}{\gamma^2}$, which are in turn equivalent to $\eta<1-\frac{\delta}{1+p_0-\frac{4}{\gamma^2}}$ and $\eta(1+\frac{4}{\gamma^2})>\delta+\frac{4}{\gamma^2}$. This is possible if and only if $(1-\frac{\delta}{1+p_0-\frac{4}{\gamma^2}})(1+\frac{4}{\gamma^2})>\delta+\frac{4}{\gamma^2}$, which produces our condition $\delta<\frac{1+p_0-\frac{4}{\gamma^2}}{2+p_0}$. Notice that the $\eta$ introduced in the proof of the upper bound is not the same as the $\eta$ introduced in the proof of the lower bound.
 \medskip
We are going to evaluate the cost term on different possible event, i.e we introduce
\begin{equation}\label{costi}
C_i(t):=\E[ \,"\, \mathbf{1}_{A_i}]
\end{equation}
where $\,"\,$ stands for the integrand inside the expectation in \eqref{cost} and the event $A_i$ ($i=1,2$) ranges respectively among the two events $\{M>\frac{\eta}{\gamma}\ln t\}$ and $\{M\leq \frac{\eta}{\gamma}\ln t\}$. 

 \medskip

Let us start with $C_1(t)$,  which can be estimated  by 
 \begin{align*}
C_1(t)
 = &\E[\frac{1}{r^{2-\gamma^2}  E_{r,v}  e^{\gamma M}}\big(\frac{I_\gamma(\infty)-I_\gamma(M)}{I_\gamma(M) I_\gamma(\infty)}\big) \mathbf{1}_{\{r^{2-\gamma^2} E_{r,v} e^{\gamma M}I_\gamma( M )>t,M>\frac{\eta}{\gamma}\ln t\}} ] \\
\leq & \E[\frac{1}{ e^{2\gamma M}} \frac{B}{I_\gamma(M)^2} F\big(\frac{t}{e^{\gamma M}I_\gamma( M )}\big)\mathbf{1}_{\{M>\frac{\eta}{\gamma}\ln t\}} ] 
\end{align*}
where we have introduced the function
\begin{equation}\label{def:F}
F(u)=\E[\frac{1}{r^{2-\gamma^2}  E_{r,v} }\mathbf{1}_{ \{r^{2-\gamma^2}  E_{r,v}>u\} }].
\end{equation}
The above expression for the cost will be analyzed according to different possible regimes of the function $F$, depending on the possible values of its argument $\frac{t}{e^{\gamma M}I_\gamma( M )}$. For this, observe that $1<2\gamma Q(1-\frac{\sqrt{2}}{Q})$, hence we can choose $0<q<2$ and $0<a<1$ such that
\begin{equation}\label{defa}
1<aq\gamma Q \quad \text{ and }\quad   0< a<1-\frac{\sqrt{2}}{Q} .
\end{equation}

Now we restrict the cost $C_1(t)$ further to the events $\{tr^{-a\gamma Q}\leq e^{\gamma M}I_\gamma( M )\}$ and then $\{tr^{-a\gamma Q}>e^{\gamma M}I_\gamma( M )\}$, producing two quantities that we respectively call $C^1_1(t)$ and $C_1^2(t)$.

Concerning $C_1^1(t)$, on the event $\{tr^{-a\gamma Q}\leq e^{\gamma M}I_\gamma( M )\}$, we can use the rough estimate $F(u)\leq r^{-2}$ (obtained by using the fact that indicator functions are bounded by $1$ in \eqref{def:F}) to deduce that  for any $0<q<2$ 
 \begin{align}
C_1^1(t) \leq & r^{-2} \E[   \frac{B}{( e^{\gamma M}I_\gamma(M) )^2}  \mathbf{1}_{\{M>\frac{\eta}{\gamma}\ln t\}}\mathbf{1}_{\{  tr^{-a\gamma Q}<e^{\gamma M}I_\gamma( M )   \}} ] \nonumber\\
=&r^{qa \gamma Q -2} t^{-q}\E[ \frac{B}{( e^{\gamma M}I_\gamma(M) )^{2-q}} \mathbf{1}_{\{M>\frac{\eta}{\gamma}\ln t\}}\mathbf{1}_{\{  tr^{-a\gamma Q}<e^{\gamma M}I_\gamma( M )   \}} ]\nonumber\\
\leq &Cr^{qa \gamma Q -2} t^{-q}\E[ \frac{1}{( e^{\gamma M} )^{2-q}} \mathbf{1}_{\{M>\frac{\eta}{\gamma}\ln t\}}  ].\label{S2}
\end{align}
Above we have used the fact that $\E[\frac{B}{I_\gamma(M)^{q-2}}|M]$ can be bounded by constant. Indeed, we use H\"older's inequality with conjugate exponents $m,m'$ and   $B\amalg M$ as well as $I_\gamma(0)\amalg M$  to get the fact that 
$$\E[ \frac{B}{I_\gamma(M)^{2-q}}|M]\leq \E[ \frac{B}{I_\gamma(0)^{2-q}}|M]\leq \E[B^m]^{1/m}\E[I_\gamma(0)^{-(2-q)m'}]^{1/m'}<+\infty$$ provided that $m$ is chosen smaller than $\frac{4}{\gamma^2}$ (recall that $I_\gamma(0)$ has negative moments of all order: see \cite{KRV}). Now we can use the explicit exponential law for $M$ to get that \eqref{S2} is less than $Cr^{qa\gamma Q -2} t^{-q-(\frac{2Q}{\gamma}-q)\eta} $. Condition \eqref{defa} imposes that  the $r$-exponent   satisfies $qa \gamma Q -2>-1$. The $t$-contribution equals $t^{-(1+\frac{4}{\gamma^2})\eta-q(1-\eta)}$, which is as expected thanks to condition \eqref{condition}.

\medskip 
Concerning $C_1^2(t)$, 
by using the Girsanov transform with the term $  e^{-\gamma N_{v,r} -\frac{\gamma^2}{2}\E[ N_{v,r}^2  ]} $, we get for some standard Gaussian random variable $Z$ and for $ur^{-a\gamma Q}>1$
\begin{align*}
F(u)=&r^{-2}\P\Big(Z>Q(-\ln r)^{1/2}(1+\frac{1}{\gamma Q}\frac{\ln u}{-\ln r})\Big) \leq r^{-2}\P\Big(Z>(1-a)Q(-\ln r)^{1/2} \Big)\leq Cr^{\frac{1}{2}(1-a)^2Q^2-2}.
\end{align*}
Then, using this estimate, $C_1^2(t)$ is seen to be less than 
 $$Cr^{\frac{1}{2}(1-a)^2Q^2-2}\E[ \frac{1}{( e^{\gamma M}I_\gamma(M) )^{2}} \mathbf{1}_{\{M>\frac{\eta}{\gamma}\ln t\}}  ]$$
 which is less than $Cr^{\frac{1}{2}(1-a)^2Q^2-2}t^{-(1+\frac{4}{\gamma^2})\eta}$ using  the explicit exponential distribution \eqref{tailofmax} of $M$ as in \eqref{S2}. Again, conditions \eqref{defa} and   \eqref{condition} ensure respectively that the $r$-exponent and $t$-exponent behave as expected. This concludes the case of $C_1(t)$.
 
 \medskip
Now we analyze $C_2(t)$
 \begin{align*}
C_2(t)\leq & \E[ \frac{1}{r^{2-\gamma^2}  E_{r,v}  e^{\gamma M}}\big(\frac{I_\gamma(\infty)-I_\gamma(M)}{I_\gamma(M) I_\gamma(\infty)}\big) \mathbf{1}_{\{r^{2-\gamma^2} E_{r,v} e^{\gamma M}I_\gamma( M )>t,M\leq \frac{\eta}{\gamma}\ln t\}} ] \\
\leq &t^{-p}\E[(r^{2-\gamma^2}  E_{r,v} )^{p-1}] \E[ ( e^{ \gamma M} I_\gamma(M))^{p-1} \mathbf{1}_{\{M\leq \frac{\eta}{\gamma}\ln t\}} ] \\
\leq &C r^{\psi(p-1)}t^{-p} \E[   I_\gamma(\infty)^{p-1}] \E[ ( e^{ \gamma M}  \mathbf{1}_{\{M\leq \frac{\eta}{\gamma}\ln t\}} ] \\
\leq &  r^{\psi(p-1)} t^{-p+\eta(p-\frac{4}{\gamma^2})},
\end{align*}
where we have used that almost surely $\frac{I_\gamma(\infty)-I_\gamma(M)}{ I_\gamma(\infty)} \leq 1$. Hence condition \eqref{condition} ensures that the cost term is $ r^{\psi(p-1)} o(t^{-\frac{4}{\gamma^2}-\delta})$.
%
%
%
To conclude it suffices to bound the term
\begin{align*}
\E[&\frac{1}{r^{2-\gamma^2}  E_{r,v}  e^{\gamma M}I_\gamma(\infty)} \mathbf{1}_{\{r^{2-\gamma^2} E_{r,v} e^{\gamma M}I_\gamma( M )>t\}} ] \\
&\leq \E[\frac{1}{r^{2-\gamma^2}  E_{r,v}  e^{\gamma M}I_\gamma(\infty)} \mathbf{1}_{\{r^{2-\gamma^2} E_{r,v} e^{\gamma M}I_\gamma( \infty )>t\}} ] \\
&=\frac{2(Q-\gamma)}{2(Q-\gamma)+\gamma}\bar{R}(\gamma)t^{-\frac{4}{\gamma^2}}+r^{\psi(p-1)}o(t^{-\frac{4}{\gamma^2}-\delta})
\end{align*}
where $o(t^{-\frac{4}{\gamma^2}-\delta})$ is uniform in $v$,$r$. The last line is obtained by independence of $E_{r,v}$, $M$, $I_\gamma(\infty)$  and the explicit tail of the exponential distribution. Our claim follows.\qed

\medskip

\noindent {\it Proof of Corollary \ref{th2}.}  In what follows, the Gaussian process $X_0$ stands for the GFF inside $D$ with vanishing mean on the unit circle, i.e.  with covariance structure given by \eqref{hatGformulafirst} whereas $G$ is an independent centered Gaussian field in $D$ such that 
$$\E[G(x)G(y)]=f(x,y),$$ with $f$ the function appearing in \eqref{defkernel}. Then we may assume that $X=X_0+G$.

Now we follow and adapt the proof of Theorem \ref{th1}. For this, we set
$$M^0_\gamma(dx)= e^{\gamma X_0(x)-\frac{\gamma^2}{2} \E[X_0(x)^2]} dx$$
in such a way that $M_\gamma(dx)=e^{\gamma G(x)-\frac{\gamma^2}{2}\E[G(x)^2]}M^0_\gamma(dx)$.
The localization trick then yields
\begin{align*}
\P(M_\gamma(\mathcal{O})>t)&=  \E[\frac{M^0_\gamma(\mathcal{O})}{M^0_\gamma(\mathcal{O})}\mathbf{1}_{\{M_\gamma(\mathcal{O})>t\}}]=\int_{\mathcal{O}}  \E[e^{\gamma X_0(v)-\frac{\gamma^2}{2}\E[X_0(v)^2]}\frac{1}{M^0_\gamma(\mathcal{O})}\mathbf{1}_{\{M_\gamma(\mathcal{O})>t\}}]\,d^2v
\end{align*}
and Girsanov's transform asserts that weighting the probability law by $e^{\gamma X_0(v)-\frac{\gamma^2}{2}\E[X_0(v)^2]}$ amounts to shifting the law of $X$ by $\gamma \ln\frac{1}{|\cdot-v|}$, hence
\begin{align}\label{ajout}
\P(M_\gamma(\mathcal{O})>t)&=  \int_\mathcal{O}\E[\frac{1}{M^0_\gamma(v,\mathcal{O}) } \mathbf{1}_{\{M_\gamma(v,\mathcal{O})>t\}}] d^2v
\end{align}
where we have set 
$$M^0_\gamma(v,\mathcal{O}):=\int_{\mathcal{O}}\frac{1}{|z-v|^{\gamma^2}}M^0_\gamma(d^2z)\quad \text{ and }\quad M_\gamma(v,\mathcal{O}):=\int_{\mathcal{O}}\frac{1}{|z-v|^{\gamma^2}}M_\gamma(d^2z).$$
Let us set $r={\rm dist}(v,\partial \mathcal{O})$ and $r'=\min(r,\epsilon)$, where $\epsilon>0$ is a regularization parameter which will be sent to $0$ in the end.
Sticking to the notations of subsection \ref{rid}, we set
$$A_0:=M^0_\gamma(v,B(v,r')^c\cap\mathcal{O}) \quad \text{ and }\quad A:=M_\gamma(v,B(v,r')^c\cap\mathcal{O}).$$
Therefore $M_\gamma(v,\mathcal{O})=M_\gamma(v,B(v,r'))+A$ and similarly for the corresponding items with index $0$. Our next step will be to use the computations already done for items with index $0$ in the proof of Theorem \ref{th1} and compare with items $M_\gamma(v,B(v,r'))$ and $A$.  The difference between these quantities involves the process $G$, which will be estimated in terms of the quantities
\begin{align}
\mathcal{S}_{r'}(v):=\sup_{z\in B(v,r')}e^{\gamma G(z)-\frac{\gamma^2}{2}\E[G(z)^2]} & & \mathcal{S}:=\sup_{z\in \mathcal{O}}e^{\gamma G(z)-\frac{\gamma^2}{2}\E[G(z)^2]} \\
\mathcal{I}_{r'}(v):=\inf_{z\in B(v,r')}e^{\gamma G(z)-\frac{\gamma^2}{2}\E[G(z)^2]} & & \mathcal{I}:=\inf_{z\in \mathcal{O}}e^{\gamma G(z)-\frac{\gamma^2}{2}\E[G(z)^2]}. \label{Gmin}
\end{align}
Then we can reproduce the argument \eqref{upperbound} for the upper bound
\begin{align}
\E[\frac{1}{M^0_\gamma(v,\mathcal{O}) } \mathbf{1}_{\{M_\gamma(v,\mathcal{O})>t\}}]\leq & \E[\frac{1}{M^0_\gamma(v,B(v,r'))+A } \mathbf{1}_{\{M_\gamma(v,B(v,r'))>t-t^{1-\eta}\}}]+\E[\frac{1}{A } \mathbf{1}_{\{A>t^{1-\eta}\}}] \nonumber\\
\leq & \E[\frac{1}{M^0_\gamma(v,B(v,r'))} \mathbf{1}_{\{M_\gamma(v,B(v,r'))>t-t^{1-\eta}\}}]+\E[A^{p-1}]t^{-(1-\eta)p} \nonumber\\
\leq & \E[\frac{1}{M^0_\gamma(v,B(v,r'))} \mathbf{1}_{\{\mathcal{S}_{r'}(v) M^0_\gamma(v,B(v,r))>t-t^{1-\eta}\}}]+\E[A^{p-1}]t^{-(1-\eta)p} \label{refAetc}
\end{align}

We have the obvious bound for $p\in (0,\tfrac{4}{\gamma^2})$
$$\E[A^p]\leq \E[\mathcal{S}^p]\E[A^p_0].$$
The assumption of $f$ being  H\"older on $\bar{\mathcal{O}}$ ensures that $ \E[\mathcal{S}^p]$ is finite by standard arguments for the supremum of Gaussian processes (see lecture 6 in \cite{Biskup} for example). Hence $A$ satisfies the same bounds as Lemma \ref{lemA} and the second term in \eqref{refAetc} is again of the form $(r')^{\psi(p-1)}o(t^{-\frac{4}{\gamma^2}-\delta})$.

The first term in \eqref{refAetc} concerning the behaviour near the singularity at $v$ is  bounded as in  subsection \ref{nearsing} which states that
$$ \E[\frac{1}{M^0_\gamma(v,B(v,r'))} \mathbf{1}_{\{ M^0_\gamma(v,B(v,r'))>t-t^{1-\eta}\}}]\leq \frac{2(Q-\gamma)}{2(Q-\gamma)+\gamma}\bar{R}(\gamma)t^{-\frac{4}{\gamma^2}} + (r')^{-\alpha}o(t^{-\frac{4}{\gamma^2}-\delta})$$ for some $\alpha\in (0,1)$.
By independence of $\mathcal{S}_{r'} (v)$ and conditioning on $\mathcal{S}_{r'} (v)$, we deduce
$$ \E[\frac{1}{M^0_\gamma(v,B(v,r'))} \mathbf{1}_{\{\mathcal{S}_{r'}(v) M^0_\gamma(v,B(v,r'))>t-t^{1-\eta}\}}]\leq   \frac{2(Q-\gamma)}{2(Q-\gamma)+\gamma}\bar{R}(\gamma)t^{-\frac{4}{\gamma^2}}\E[(\mathcal{S}_{r'} (v))^{\frac{4}{\gamma^2}}] + (r')^{-\alpha}o(t^{-\frac{4}{\gamma^2}-\delta}).$$
Integrating this relation over $\mathcal{O}$ we deduce from \eqref{ajout} that
$$\limsup_{t\to\infty}t^{\frac{4}{\gamma^2}}\P(M_\gamma(\mathcal{O})>t)\leq \left(\int_{\mathcal{O}}\E[(\mathcal{S}_{r'} (v))^{\frac{4}{\gamma^2}}] \,dv\right)\frac{2(Q-\gamma)}{2(Q-\gamma)+\gamma}\bar{R}(\gamma).$$
This bound is valid for arbitrary $\epsilon$ (recall that $r'=\min(r,\epsilon))$. So we want to let $\epsilon\to 0$. An easy application of the dominated convergence theorem ensures that
$$\lim_{\epsilon\to 0}\int_{\mathcal{O}}\E[(\mathcal{S}_{r'} (v))^{\frac{4}{\gamma^2}}] \,dv=\int_{\mathcal{O}}\E[(e^{\gamma G(v)-\frac{\gamma^2}{2}\E[G(v)^2]})^{\frac{4}{\gamma^2}}] \,dv=\int_{\mathcal{O}}e^{\frac{4}{\gamma}(Q-\gamma)f(v,v)}\,dv.$$
Hence 
$$\limsup_{t\to\infty}t^{\frac{4}{\gamma^2}}\P(M_\gamma(\mathcal{O})>t)\leq \left(\int_{\mathcal{O}}e^{\frac{4}{\gamma}(Q-\gamma)f(v,v)}\,dv\right)\frac{2(Q-\gamma)}{2(Q-\gamma)+\gamma}\bar{R}(\gamma).$$
This shows the upper bound. The lower bound is established in the same way (using \eqref{Gmin}).\qed

\section{Extensions to other cases}\label{beyond}
In this section, we explain how to generalize our results to the other cases. Though the claims of this section cannot be taken for granted (due to the fact that it sweeps under the rug potential technical difficulties), it provides nonetheless arguments which we believe are convincing to tackle the other cases.

\subsection{The 2d case}
Consider now the general case in 2d of a log-correlated kernel of the type \eqref{defkernel}. Along the same lines as for the proof of Theorem \eqref{th1} the  localization trick  allows us to trade the study of the tail of the random variable $M_\gamma(\mathcal{O})$ for the study of the tail of the singular integral (for some $r>0$)
\begin{equation}\label{tailgeneral2d}
M_\gamma(v,B(v,r)):= \int_{B(v,r)} e^{\gamma^2 f(v,z)}  \frac{e^{    \gamma X(z)-\frac{\gamma^2}{2} \E[X(z)^2]   }}{|z-v|^{\gamma^2}} d^2z
\end{equation}
for each point $v \in \mathcal{O}$  ($f$ appears in \eqref{defkernel}). The full tail of $M_\gamma(\mathcal{O})$ is obtained by integration of the term 
$$ \E[\frac{1}{M_\gamma(v,B(v,r))} \mathbf{1}_{\{M_\gamma(v,B(v,r))>t\}}]$$
 with respect to the measure $d^2v$. The advantage of this localization trick is that the  tail of \eqref{tailgeneral2d} is completely concentrated around $v$, hence is determined by a local analysis which is more transparent than a direct  study of $M_\gamma(\mathcal{O}) $. 
For $r$ small, the variable $X$ has roughly covariance of the form $\ln \frac{1}{|x-y|}+ f(v,v)$ and therefore $X \approx \bar{X}+N$ (in law) where $\bar{X}$ has $\ln \frac{1}{|x-y|}$ covariance and $N$ is a centered Gaussian independent from everything and with variance $f(v,v)$. Combining this decomposition with our estimates for $\bar{X}$, one expects that 
\begin{equation}\label{tailgeneral2dresult}
t^{\frac{4}{\gamma^2}} \E[\frac{1}{M_\gamma(v,B(v,r)) } \mathbf{1}_{\{M_\gamma(v,B(v,r))>t\}}] \underset{t \to \infty}{\rightarrow} \frac{\frac{2}{\gamma}(Q-\gamma)}{\frac{2}{\gamma}(Q-\gamma)+1} \bar{R}(\gamma)    e^{ \frac{4}{\gamma}(Q-\gamma) f(v,v) }. 
\end{equation}
 
 Gathering the above considerations, we expect the following generalization of Theorem \ref{th1} for kernels of type \eqref{defkernel}
 \begin{equation}
 \P(M_\gamma(\mathcal{O})>t)=   \left (  \int_\mathcal{O}  e^{ \frac{4}{\gamma}(Q-\gamma) f(v,v) } d^2v \right )   \frac{\frac{2}{\gamma}(Q-\gamma)}{\frac{2}{\gamma}(Q-\gamma)+1} \bar{R}(\gamma)  \frac{1}{t^{\frac{4}{\gamma^2}}} +o(t^{-\frac{4}{\gamma^2}}).
 \end{equation}

We leave open the determination of bounds on the $o(t^{-\frac{4}{\gamma^2}})$ term.

\subsection{The other dimensions} 
In higher dimensions, we expect the method to work as well by decomposing the log-correlated field into a radial part (Brownian motion) and an independent radial part around each localization point $v$. The constant in the expansion will then be given by explicit terms times  the $d$-dimensional analog $\bar{R}_{{\rm d}}(\gamma)$ of the reflection coefficient $\bar{R}(\gamma)$ defined by \eqref{firstdefunitR}. The question is then to know if we can compute explicitly this expectation depending on $d$; presently, getting explicit formulas in dimension $d \geq 3$ seems out of reach since one can not rely on the powerful framework of $2d$ conformal field theory for the GFF.

Like in dimension $2$, we also have an explicit expression for $\bar{R}_{{\rm 1}}(\gamma)$ in dimension $1$. Indeed, one has the following expression for $\bar{R}_1 (\gamma)$ (the so-called boundary unit volume reflection coefficient in the terminology of Liouville field theory):
\begin{equation*}
\bar{R}_1(\gamma)=  \E[ \left ( \int_{-\infty}^\infty  e^{   \gamma \mathcal{B}_s^{1,\gamma}}  e^{\gamma Y(s)-\frac{\gamma^2 E[Y(s)^2]}{2}}   Z_s ds  \right )^{\frac{2}{\gamma} (Q_1-\gamma)}]
\end{equation*}
where $\mathcal{B}_s^{1,\gamma} $ is defined like $\mathcal{B}_s^\gamma $ in \eqref{BMneg} with $Q$ replaced by $Q_1=\frac{\gamma}{2}+\frac{1}{\gamma}$ and $Y_s$ is the restriction to the real line of the centered Gaussian field with covariance \eqref{covlateral}. This yields the asymptotic
$$\P(M_\gamma(\mathcal{O})>t)= \left (  \int_\mathcal{O}  e^{ \frac{2}{\gamma}(Q_1-\gamma) f(v,v) } d^2v \right )(1-\tfrac{\gamma^2}{2})  \frac{\bar{R}_{{\rm 1}}(\gamma)}{t^{\frac{2}{\gamma^2}}} +o(t^{-\frac{2}{\gamma^2}}).$$
The recent integrability results of R\'emy for GMC on the circle \cite{remy} allows us to compute explicitly the tail of $M_\gamma(\mathcal{O})$ in the case when $X$ is the circular logarithmic noise and $\mathcal{O}=(0,2\pi)$; as a matter of fact, the result of R\'emy is much more precise since it gives the precise density of the total mass of GMC on the circle (this density was conjectured in the physics literature in 2008 by Fyodorov-Bouchaud \cite{fybu}; see also a similar conjecture in \cite{FLeR} for the case of the unit interval). This leads to the  following explicit expression for
$\bar{R}_1(\gamma)$
\begin{equation*}
\bar{R}_1(\gamma)= \frac{(2\pi)^{ \frac{2}{\gamma}(Q_1-\gamma)}}{(1-\tfrac{\gamma^2}{2})\Gamma(1-\tfrac{\gamma^2}{2})^{\frac{2}{\gamma^2}}}.
\end{equation*}


\end{document}